\documentclass{article}
\usepackage{amstex,psfig}

\newtheorem{theorem}{Theorem}
\bigskip
\begin{document}

\title{Mathematical Modeling of Boson-Fermion Stars in the Generalized
     Scalar-Tensor Theories of Gravity\thanks{
This research was supported by the Bulgarian Ministry of Education, Science
and Technologies under the Grants NoNo MM-602/96, F610/99 and by the Sofia
University Research Fund, Contr. No 245/99.}}
\author{T.~L. Boyadjiev \\
{\small \textit{Faculty of Mathematics and Computer Science
}}\\
[-1.mm] {\small \textit{University of Sofia, 1164 Sofia, Bulgaria}}\\
[-1.mm] {\small \textit{e-mail: todorlb@@fmi-uni.sofia.bg}}\\
\\
M.~D. Todorov\\
{\small \textit{Faculty of Applied Mathematics and Computer Science}}\\
[-1.mm] {\small \textit{Technical University of Sofia, 1756 Sofia,
Bulgaria}}\\
[-1.mm] {\small \textit{e-mail: mtod@@vmei.acad.bg}}\\
\\
P.~P. Fiziev\\
{\small \textit{
Faculty of Physics, University of Sofia,
1164 Sofia, Bulgaria}}\\
[-1.mm] {\small \textit{e-mail: fiziev@@phys.uni-sofia.bg}}\\
\\
S.~S. Yazadjiev\\
{\small \textit{
Faculty of Physics, University of Sofia, 1164 Sofia,
Bulgaria}}\\
[-1.mm] {\small \textit{e-mail: yazad@@phys.uni-sofia.bg}}}
\date{}
\maketitle

\begin{abstract}
A model of static boson-fermion star with spherical
symmetry based on the scalar-tensor theory of gravity with massive
dilaton field is investigated numerically.

Since the radius of star is \textit{a priori} an unknown quantity,
the corresponding boundary value problem (BVP) is treated as a
nonlinear spectral problem with a free internal boundary. The
Continuous Analogue of Newton Method (CANM) for solving this
problem is applied.

Information about basic geometric functions and the functions
describing the matter fields, which build the star is obtained.
In a physical point of view the main result is that the structure
and properties of the star in presence of massive dilaton field
depend essentially both of its fermionic and bosonic components.


\end{abstract}

\bigskip

\section{Introduction}
The most natural and promising generalizations of general
relativity are the scalar-tensor theories of gravity \cite{BD} --
\cite{DEF1}. In these theories gravity is mediated not only by a
tensor field (the metric of space-time) but also by a scalar
field (the dilaton). The scalar-tensor theories of gravity
contain arbitrary functions of the scalar field that determine
the gravitational ``constant" as a dynamical variable and the
strength of the coupling between the scalar field and matter. It
should be stressed that specific scalar-tensor theories of
gravity arise naturally as a low energy limit of the string
theory \cite{GSW} -- \cite{MAHS}, which is the most promising
modern model of the unification of all fundamental physical
interactions.

If the string theory and its low energy limit are relevant to the
real world, then the dilaton must be massive \cite{DP}.
Unfortunately, our current understanding of how the dilaton
acquires mass is primitive and it is tied to our lack of
understanding of supersymmetry breaking. At present, we do not
have a model of how the dilaton mass is generated in the string
theory. Besides the mass term for dilaton field we may consider
the general case of arbitrary dilaton potential, describing its
nonlinear self-interaction.

From physical point of view, it is important to know how the
dilaton mass and, in general, the dilaton potential influence the
structure and stability of compact objects such as neutron stars,
boson stars, and mixed fermion-boson stars.

It is known that the predictions of scalar-tensor theories of
gravity with massless dilaton may differ drastically from these
of general relativity. For example, the phenomenon of
``spontaneous scalarization" was discovered recently \cite{DEF2},
\cite{SSN} as a non-perturbative strong field effect in a massive
neutron star. The existence of this effect poses some important
physical questions \cite{DEF3}. That is why it is natural to ask
whether or not the ``spontaneous scalarization" will occur when
the dilaton is massive. In recent years, the boson stars in
scalar-tensor theories of gravity with massless dilaton have been
widely studied both analytically and numerically (see for example
\cite{GJ} - \cite{CS}). The study of boson stars in the case of
massive dilaton is physically interesting and may be important
for the understanding of their formation in the early universe.

The investigation of the compact objects in the generalized
scalar-tensor theories of gravity helps us understand them
better. On the other hand, the investigation of matter in extreme
conditions like these in the neutron stars may demonstrate new
phenomena and new features of specific scalar-tensor theories of
gravity, originating from the low energy limit of the string
theory. Thus, at first time we may be able to reach theoretical
indications of physical manifestation of the string theory in the
real world \cite{BFY}.

In the present paper we develop a direct numerical method for
solving the equations of the general scalar-tensor theories of
gravity including a dilaton potential term for the general case of
mixed boson-fermion star.

The physical motivation for considering mixed boson-fermion stars
is connected with the fact that many of the present-day existing
stars are of primordial origin being formed from an original gas
of fermions and bosons in the early universe. That is why it
should be expected that they are a mixture of both fermions and
bosons in different proportions. The study of such mixed objects
is a new interesting problem, whose investigation was started in
reference \cite{HLM}. There exist different candidates for boson
fields in stars such as Higgs field of Standard model, or axion
field being a pseudoscalar partner of dilaton in the superstring
theory. They are unavoidable part of modern physics, nevertheless
up to now we have no experimental evidence for their existence.
Taking into account that according to the modern understanding of
the initial state of universe a significant amount of these
fields must have been present during the Big Bang phase, one has
to expect some part of these fields to be present in the stars of
primordial origin. The study of new observable effects of boson
fields in such mixed stars may give new ways for discovery of the
existence of the above hypothetical fields, which at present are
the most intriguing new objects in modern physics.

In the Einstein frame the field equations in the presence of
fermion and boson matter are:
\begin{eqnarray}
G_{i }^{j }=\kappa _{*}\left( \stackrel{B}{T_{i }^{j }}+\stackrel{F}{%
T_{i }^{j }}\right) +2\partial _{i}\varphi \partial ^{j }\varphi
\!\!\! &&
-\partial ^{l }\varphi \partial _{l }\varphi \delta _{i }^{j }+%
\frac{1}{2}U(\varphi )\delta _{i}^{j},  \nonumber   \\ \nabla _{i
}\nabla ^{i }\varphi +{\frac{1}{4}}U^{\prime }(\varphi )\!\!\!&&
=-
\frac{\kappa _{*}}{2} \> \alpha
\varphi )\left( \stackrel{B}{T}+\stackrel{F}{T}%
\right), \nonumber\\ [-2mm] \label{SFE}\\[-2mm] \nabla _{i }\nabla
^{i }\Psi +2\alpha (\varphi )\partial ^{l
}\varphi \partial _{l }\Psi \!\!\! && =-2A^{2}(\varphi ){\frac{\partial
{\tilde{%
W}}}{\partial \Psi ^{+}}}, \nonumber \\ \nabla _{i }\nabla ^{i
}\Psi ^{+}+2\alpha (\varphi )\partial ^{l
}\varphi \partial _{l }\Psi ^{+} \!\!\!&& =-2A^{2}(\varphi ){\frac{\partial
{%
\tilde{W}}}{\partial \Psi ^{+}}},  \nonumber
\end{eqnarray}
\noindent where $\nabla _{i }$ is the Levi-Civita connection with
respect to the metric $g_{i j },(i =0,...,3;j =0,...,3)$. The
constant $\kappa _{*}$ is given by $\kappa _{*}=8\pi G_{*}$, where
$G_{*}$ is the bare Newtonian gravitational constant. The physical
gravitational ``constant" is $G_{*}A^{2}(\varphi )$, where
$A(\varphi )$ is a function of the dilaton field $\varphi $
depending on the concrete scalar-tensor theory of gravity. For
example, in the framework of the Brans-Dicke model we have
$A(\varphi )=\exp (\frac{\varphi }{\sqrt{2\omega_{BD}+3}})$,
where $\omega_{BD}$ is a parameter.

The dilaton potential $U(\varphi )$ can be written in the form
$U(\varphi )=m_{D}^{2}V(\varphi )$, where $m_{D}$ is the dilaton
mass and $V(\varphi )$ is a dimensionless model function of
$\varphi $.

The complex scalar field $\Psi $ describes the bosonic matter,
while $\Psi ^{+}$ is its complex conjugated function. The
quantity $W(\Psi ^{+}\Psi )$ is the potential of boson field,
which can be chosen in the following form:
\[
{\tilde{W}}(\Psi ^{+}\Psi )=-\frac{m_{B}^{2}}{2}\Psi ^{+}\,\Psi -\frac{1}{4%
}{\tilde{\Lambda}\,}(\Psi ^{+}\Psi )^{2},
\]
\noindent where $\tilde{\Lambda}$ is a parameter.

The scalar function $\alpha (\varphi )=\frac{d}{d\varphi }\left[
\ln A(\varphi )\right]$ determines the strength of the coupling
between the dilaton field $\varphi $ and matter.

The quantities $\stackrel{B}{T}$ and $\stackrel{F}{T}$ are
correspondingly the trace of the energy-momentum tensor of the
fermionic matter $\stackrel{F}{T_{i }^{j }}$  and the bosonic
matter $\stackrel{B}{T_{i }^{j}}$. We note that in the present
article we consider the fermionic matter only in macroscopic
approximation, {\it i.e.}, after averaging quantum fluctuations
of the corresponding fermion fields. Thus, we actually consider
standard classical relativistic matter.

The explicit forms of the mentioned tensors are correspondingly:
\begin{eqnarray}
\stackrel{B}{T_{i}^{j}}=&&\!\!\! \frac{1}{2}A^{2}(\varphi )\left(
\partial _{i}\Psi ^{+}\partial ^{j}\Psi +\partial _{i
}\Psi \partial ^{j}\Psi ^{+}\right)  \nonumber \\ &&\!\!\! \qquad
-\frac{1}{2}A^{2}(\varphi )\left[ \partial _{l}\Psi ^{+}\partial
^{l}\Psi -2A^{2}(\varphi ){\tilde{W}}(\Psi ^{+}\Psi )\right]
\delta _{i}^{j}\>, \label{eq:2a} \\ \stackrel{F}{T_{i}^{j
}}=&&\!\!\! \left( \varepsilon +p\right) u_{i}\,u^{j}-p\,\delta
_{i}^{j}\>. \label{eq:2b}
\end{eqnarray}
\noindent Here, the energy density and the pressure of the
fermionic fluid in the Einstein frame are $\varepsilon
=A^{4}(\varphi ){\tilde{\varepsilon}}$ and $p=A^{4}(\varphi
)\,\tilde{p}$, where ${\tilde{\varepsilon}}$ and $\tilde{p}$ are
the physical energy density and pressure. Instead of giving the
equation of state of the fermionic matter in the form
$\tilde{p}=\tilde{p}({\tilde{\varepsilon}})$, it is more
convenient to write it in a parametric form:
\begin{equation}
{\tilde{\varepsilon}}={\tilde{\varepsilon}_{0}}g(\mu )\,\,\,\,\,\,\,\,\,\,%
\tilde{p}={\tilde{\varepsilon}_{0}\,}f(\mu ),  \label{eq:2c}
\end{equation}

\noindent where ${\tilde{\varepsilon}_{0}}$ is a properly chosen
dimensional constant, $\mu $ is the dimensionless Fermi momentum,
and $f(\mu)$ and $g(\mu)$ are given functions (see below).

The physical four-velocity of the fermionic fluid is denoted by
$u_{i}$.

The field equations together with the Bianchi identities lead to
the local conservation law of the energy-momentum of matter:
\begin{equation}
\nabla _{j}\stackrel{F}{T_{i}^{j}}=\alpha (\varphi )\stackrel{F}{T}%
\partial _{i}\varphi \,.  \label{TBIAN}
\end{equation}

From now on, we will take into consideration a static and
spherically symmetric mixed boson-fermion star in asymptotic flat
space-time. This means that the metric $g_{ij}$ has the form:
\begin{equation}
ds^{2}=e^{\nu (r)}dt^{2}-e^{\lambda (r)}dr^{2}-r^{2}\left(
d\theta ^{2}+\sin ^{2}\theta \,d\phi ^{2}\right), \label{eq:gmn}
\end{equation}
where $r,\theta ,\phi $ are usual spherical coordinates.

The field configuration is static when the boson field $\Psi $
satisfies the condition:
\[
\Psi ={\tilde{\sigma}(}r{)\,}e^{i\omega t}.
\]
Here, $\omega $ is a real number and ${\tilde{\sigma}(}r)$ is a
real function.

Taking into account the above-stated assumption, the system of
the field equations is reduced to a system of ordinary
differential equations (ODEs). Before writing the system
explicitly, we are going to introduce a rescaled (dimensionless)
radial coordinate by $r\to m_{B}r$, $r\in [0,\infty)$, where
$m_{B}$ is the mass of the bosons (a prime will denote the
differentiation with respect to the dimensionless radial
coordinate $r$).

We also define the following dimensionless quantities by:
\[
\Omega ={\frac{\omega }{m_{B}}},\quad \sigma =\sqrt{\kappa }_{*}\,{\tilde{%
\sigma}},\quad \Lambda ={\frac{{\tilde{\Lambda}}}{\kappa _{*}{m^{2}}_{B}}}%
,\quad \gamma ={\frac{m_{D}}{m_{B}}}.
\]
The components of the energy-momentum tensors of the fermionic and
bosonic matter, written by the dimensionless quantities, are
correspondingly:
\begin{eqnarray}
\stackrel{\mathit{F}}{T_{0}^{0}}\!\!\!&&=bA^{4}(\varphi )\,g(\mu ),\quad
\stackrel{%
\mathit{F}}{T_{1}^{1}}\>=\>\stackrel{\mathit{F}}{T_{2}^{2}}\>=\>-bA^{4}(
\varphi
)\,f(\mu ),  \label{eq:3a} \\
\stackrel{\mathit{B}}{T_{0}^{0}}\!\!\!&&=\frac{1}{2}\Omega
^{2}A^{2}(\varphi )\,e^{-\nu }\sigma
^{2}(r)+\frac{1}{2}A^{2}(\varphi )\,e^{-\lambda }\sigma
^{\prime 2}-A^{4}(\varphi )W(\sigma ^{2}),  \label{eq:3b} \\
\stackrel{\mathit{B}}{T_{1}^{1}}\!\!\!&&=-\frac{1}{2}\Omega
^{2}A^{2}(\varphi )\,e^{-\nu }\sigma
^{2}(r)-\frac{1}{2}A^{2}(\varphi )\,e^{-\lambda }\sigma
^{\prime 2}-A^{4}(\varphi )W(\sigma ^{2}),  \label{eq:3c} \\
\stackrel{\mathit{B}}{T_{2}^{2}}\!\!\!&&=-\frac{1}{2}\Omega
^{2}A^{2}(\varphi )\,e^{-\nu }\sigma
^{2}(r)+\frac{1}{2}A^{2}(\varphi )\,e^{-\lambda }\sigma ^{\prime
2}-A^{4}(\varphi )W(\sigma ^{2}).  \label{eq:3d}
\end{eqnarray}
The parameter
$b={\kappa_{*}{\tilde{\varepsilon}_{0}}}/{m_{B}^{2}} $ describes
the relation between the Compton length of dilaton and the usual
radius of neutron star in general relativity.

It is necessary to note that two physically interesting borderline
cases of pure bosonic and pure fermionic stars are formally
contained in the above general system (\ref{SFE}). For example,
the model of pure bosonic stars can be obtained from (\ref{SFE})
by letting the tensor  $\stackrel{\mathit{F}}{T_{i}^{j}}$ to be
zero. While the pure fermionic stars correspond to the field
$\Psi \equiv 0$. The case of pure bosonic stars in the
scalar-tensor theories of gravity with a massive dilaton has
already been discussed in our recent paper \cite{FYBT}. In the
present paper we consider the mixed boson-fermion stars.

\section{Formulation of the Problem}
Under the physical assumptions we have made, the field equations
(\ref{SFE}) can be reduced to a system of ODEs. From mathematical
point of view it is more convenient all ODEs to be of second
order. That is why  we first solve the Einstein equation $G^1_1$
for $e^{\lambda}$:
\[
e^{\lambda }=\frac{1+r\nu ^{\prime }-r^{2}{\varphi
^{\prime }}^{2}-\frac{1}{%
2}A^{2}(\varphi )r^{2}{\sigma ^{\prime }}^{2}}{1-r^{2}\left[ \stackrel{%
\mathit{F}}{T_{1}^{1}}+\frac{1}{2}\gamma ^{2}V(\varphi
)-\frac{1}{2}\Omega ^{2}A^{2}(\varphi )e^{-\nu }\sigma
^{2}-A^{4}(\varphi )W(\sigma ^{2})\right] },
\]
\noindent as a function of the quantities $\nu(x)$,
$\nu^\prime(x)$, $\sigma(x)$, $\sigma^\prime(x)$, $\varphi(x)$,
$\varphi^\prime(x)$, and the spectral parameter $\Omega$, and then
substitute the above expression in the other Einstein equations.
In this way, in terms of the dimensionless quantities, the system
of the field equations (\ref{SFE}) is reduced to the following
system of ODEs:
\begin{eqnarray}
&&\nu^{\prime \prime} + \frac{\nu^{\prime}}{r} =
\Bigl[-\frac{\nu^{\prime}}{r}
+ \bigl(\stackrel{\mathit{F}}{T_0^0} - \stackrel{\mathit{F}}{T_0^0} - 2
\stackrel{\mathit{F}}{T_2^2} + \stackrel{\mathit{B}}{T_0^0} - \stackrel{%
\mathit{B}}{T_1^1} - 2\stackrel{\mathit{B}}{T_2^2} \bigr)        \nonumber
\\
&& \qquad\qquad\qquad - \gamma^2 V(\varphi) + \frac{\nu^\prime r}{2}
bigl( %
T_1 + \gamma^2 V(\varphi) \bigr) \Bigr] e^\lambda,  \label{eq:1} \\
&&\varphi^{\prime \prime} + \frac{\varphi^{\prime}}{r} = \!\! \Bigl[-\frac{%
\varphi^{\prime}}{r} +
\frac{\alpha(\varphi)}{2}\!\!\left(\stackrel{\mathit{F}}{T%
} + \stackrel{\mathit{B}}{T} \right) \nonumber \\
&&\qquad\qquad\qquad + \frac{1}{4} \gamma^2 V^\prime
(\varphi) + \frac{\varphi^\prime r}{2} \bigl( T_1 + \gamma^2 V(\varphi) %
\bigr) \Bigr] e^\lambda,  \label{eq:2} \\
&&\sigma^{\prime \prime} + \frac{\sigma^{\prime}}{r} = - 2
\alpha(\varphi) \varphi^\prime \sigma^\prime +
\biggl[-\frac{\sigma^{\prime}}{r} - 2 A^2(\varphi)
W^\prime(\sigma^2) \sigma \nonumber \\ && \qquad\qquad\qquad - \Omega^2
e^{-\nu} \sigma + \frac{%
\sigma^\prime r}{2} \bigl(T_1 + \gamma^2 V(\varphi) \bigr)
\biggr] e^\lambda. \label{eq:3}
\end{eqnarray}
In the above equations, the potential of the bosonic matter $W$
has the form:
\[
W(\sigma^2) = -\frac{1}{2} \left(\sigma^2 + \frac{1}{2}\Lambda
\sigma^4 \right),
\]
and we suppose that $W^{\prime }(\sigma^2) \equiv
\frac{dW}{d\left( \sigma^2 \right)}$. Similarly, we set
$V^\prime(\varphi) \equiv \frac{dV}{d \varphi}$.

The quantity $T_{1}$ depends on the components of the
energy-momentum tensors of the fermionic and bosonic matter
(\ref{eq:3a})-(\ref{eq:3c}):
\[
T_{1}= \>
\stackrel{\mathit{F}}{T_{0}^{0}}+\stackrel{\mathit{F}}{T_{1}^{1}}+%
\stackrel{\mathit{B}}{T_{0}^{0}}+\stackrel{\mathit{B}}{T_{1}^{1}}\newline.
\]
The quantities $\stackrel{\mathit{B}}{T}$\ \ and
$\stackrel{\mathit{F}}{T}$\ \ represent the traces of the these
tensors, and are defined by the formulae:
\[
\stackrel{\mathit{B}}{T}=-\Omega ^{2}A^{2}(\varphi )\;e^{-\nu
}\sigma ^{2}(r)+A^{2}(\varphi )\;e^{-\lambda }\sigma ^{\prime
2}-4A^{4}(\varphi )W(\sigma ^{2}),
\]
\[
\stackrel{\mathit{F}}{T}=b\,A^{4}(\varphi )\left[ g(\mu )-3f(\mu
)\right].
\]

Correspondingly, the conservation law (\ref{TBIAN}) can be
expressed as:
\begin{equation}
\mu^\prime = - \frac{g(\mu) + f(\mu)}{f^\prime(\mu)} \left[\frac{\nu^\prime%
}{2} + \alpha(\varphi) \varphi^\prime\right].        \label{eq:4}
\end{equation}
The fermionic matter functions $f(\mu )$ and $g(\mu )$, entered
in the above relations, have the form:
\begin{eqnarray}
f(\mu ) \!\!\!&&=\frac{1}{8}\left[ (2\mu -3)\sqrt{\mu +\mu ^{2}}+3\ln
\left( \sqrt{\mu }+\sqrt{1+\mu  }\right) \right], \label{eq:a1} \\
\;g(\mu ) \!\!\!&& = \frac{1}{8}\left[ (6\mu +3)\sqrt{\mu +\mu
^{2}}-3\ln \left( \sqrt{\mu }+\sqrt{1+\mu }\right) \right].
\label{eq:a2}
\end{eqnarray}

Let us now complete the problem by adding proper boundary
conditions (BCs) to the system of differential equations
(\ref{eq:1})-(\ref{eq:4}).

The asymptotic flatness means that the function $\nu (r)\to 0$
when $r \to \infty $. On the other hand, the nonsingularity
condition at the center of the star requires the derivative $\nu
^{\prime }(0)=0$. The same condition in relation to the dilaton
field $\varphi (r)$ implies that the derivative $\varphi ^{\prime
}(0)=0$. At the same time, the function $\varphi (r)$ at the
asymptotic infinity ($r\to \infty $) must be $\varphi _{\infty }
= 0$ as it is required by the asymptotic flatness. The
nonsingularity of the bosonic density $\sigma (r)$ at the center
of the star requires the derivative $\sigma ^{\prime }(0)=0$. We
need finite mass for the star, which implies $\sigma (r) \to 0$
when $r \to \infty $ . In addition, the central value $\sigma
_c=\sigma (0)$ must be given.
Concerning the fermionic fluid, we have to give the central density
${\tilde{%
\varepsilon}}_c={\tilde{\varepsilon}}(0)$ or, equivalently, the
central value $\mu _c=\mu (0)$.

It should be noted that for the physically relevant equation of
state of the fermionic matter there must be a point $r=R_s <
\infty$, where the pressure of the fermionic matter vanishes,
\textit{i.e.}, $R_s$ is the radius of the fermionic part of the
star.

As a conclusion, from the above-mentioned physical assumptions, we
can formulate the following linear boundary conditions (BCs) for
the quantities under consideration:
\begin{eqnarray}
&& \nu ^{\prime }(0)=0,\qquad \nu (\infty )=0;        \label{eq:5} \\
&& \varphi ^{\prime }(0)=0,\qquad \varphi (\infty )=0;  \label{eq:6} \\
&& \sigma ^{\prime }(0)=0,\qquad \sigma (\infty )=0;          \label{eq:7}
\\
&& \mu (0)=\mu _c\>.        \label{eq:8}
\end{eqnarray}
\noindent Here, we denote $(\cdot) (\infty) \stackrel{\rm def}{=}
\> \lim_{r \to \infty} (\cdot) (r)$.

Apart from the unknown functions $\nu (r)$, $\sigma (r)$, $\varphi (r)$, and
$%
\mu (r)$, the equations (\ref{eq:1})-(\ref{eq:4}) also include two
unknown real parameters, $R_{s}>0$ and $\Omega $. However, the
seven BCs (\ref{eq:5})-(\ref{eq:8}) are insufficient for their
computation. In order to determine these parameters, we have to
use additional conditions. In other words, the problem may be
considered as a nonlinear eigenvalue problem, where $R_{s}$ and
$\Omega $ are considered as ``eigenvalues". For this purpose,
further on we use two physically-clear additional conditions.

The first one, given by the relation:
\begin{equation}  \label{eq:4a}
\sigma(0) = \sigma_c
\end{equation}
determines the density $\sigma_c \geq 0$ of the bosonic matter in
the star's center. The second one:
\begin{equation}
 \mu(R_s) =0, \quad 0<R_{s}<\infty  \label{eq:4b}
\end{equation}
\noindent describes the condition that the density of the
fermionic matter must vanish at the radius of the star.

Finally, we note that all the functions $\nu (r)$, $\sigma (r)$,
and $\varphi(r)$ are defined in the whole real half-line $r\in
[0,\infty)$. It is easy to see that these functions are smooth in
this interval including the point $x=1$. Whereas, the fermionic
density $\mu(r)$ is defined and smooth only inside the star,
\textit{i.e.}, $r\in [0,R_s]$.

\section{Method of Solution}
\smallskip For solving the above formulated nonlinear eigenvalue problem the
Continuous Analogue of Newton Method (CANM) (see \cite{pon72} -
\cite{zhanlav}, and comprehensive surveys \cite{jmp},
\cite{puz99}) is applied. For convenience, a brief description of
CANM can be found in the attached Appendix.

The presence of the \textit{a priori} unknown quantity $R_{s}$,
however, is an obstacle for the direct use of CANM - the problem
is the unknown internal boundary $ R_{s}$. In order to overcome
this obstacle, we introduce a new scaled coordinate $x=r/R_{s}$.
As a result, the physical domain $r\in [0,\infty )$ renders to the
domain $x\in [0,\infty )$, and the star's radius $r=R_{s}$ maps
into the fixed point $x=1$. Then the BC (\ref{eq:4b}) for $\mu (x)
$ becomes
\begin{equation}
\mu (1)=0.        \label{eq:4c}
\end{equation}

Let $x_1$ and $x_2$ be two arbitrary points in the internal
domain $[0, 1]$.  We note that for the arbitrary functions $f(\mu
),g(\mu )$, and $\alpha (\varphi ) $\ the equation (\ref{eq:4})
has a first integral, which can be presented as:
\[
\int\limits_{\mu _{1}}^{\mu _{2}}\frac{f^{\prime }(\mu )}{f(\mu )+g(\mu )}%
\;d\mu +\frac{1}{2}\left( \nu _{2}-\nu _{1}\right) +\ln \frac{A(\varphi
_{2})%
}{A(\varphi _{1})}=0,
\]
where $\nu _{1}$, $\nu _{2}$, $\varphi _{1}$, $\varphi _{2}$,
$\mu _{1}$, $\mu _{2}$\ stand for the functions $\nu (x)$,
$\varphi (x)$, $\mu (x)$ at the points $x_1$ and $x_2$,
respectively. Thus, for the model of the fermionic matter
described by the conditions (\ref{eq:a1}), (\ref{eq:a2}) we
simply get the following algebraic equation:
\begin{equation}
\ln \,\left[ \frac{(1+\mu _{2})A^{2}(\varphi _{2})}{(1+\mu
_{1})A^{2}(\varphi _{1})}\right] +\nu _{2}-\nu _{1}=0.
\label{eq:i1}
\end{equation}

For convenience, we introduce the vector $\mathbf{y}(x)$ $=\{\nu
(x),$ $\varphi (x),$ $\sigma (x)\}.$ Then the first three
equations (\ref{eq:1}) - (\ref{eq:3}) of the problem and the
corresponding BCs (\ref{eq:5}) - (\ref{eq:7}) can be rewritten as
follows:
\begin{equation}
-x\mathbf{y}^{\prime \prime }-\mathbf{y}^{\prime }+\mathbf{F}=0,
\label{eq:9}
\end{equation}
\begin{equation}\label{bc9}
  \mathbf{y}^\prime(0)=0, \quad \mathbf{y}(\infty) =0,
\end{equation}
\noindent where $\mathbf{F}=\mathbf{F}(x,\mathbf{y},
\mathbf{y}^{\prime}, \mu, R_s, \Omega )$ is 3D vector consisting
of the right-hand sides (RHSs) of the equations
(\ref{eq:1})-(\ref{eq:3}) multiplied by $R_{s}^{2}x$. The
differentiation with respect to the new independent variable $x$
is denoted by $(.)^{\prime }$. In the linear case, the advantages
of such representation of the radial operator are discussed in
\cite{vidar}.

Following CANM, we introduce a ``time-like" parameter $t\in
[0,\infty )$ and assume the unknown quantities depend on $t$ as
well: $\mathbf{y}=\mathbf{y}(x,t),\> R_{s}=R_{s}(t),\> \Omega
=\Omega (t).$ Let us suppose that the function $\mu = \mu (x)$ is
known (see below). Then the CANM equations \cite{puz99}\
corresponding to (\ref{eq:9}) and (\ref{bc9}) become:
\begin{eqnarray}
-x\mathbf{z}^{\prime \prime }\!\! &+&\!\! \left(\frac{\partial
\mathbf{F}}{\partial \mathbf{y}^{\prime }} - \mathbf{E}\right)
\mathbf{z}^{\prime } +
\frac{\partial \mathbf{F}%
}{\partial \mathbf{y}}\mathbf{z}+\left( \frac{2}{R_{s}}\mathbf{F}+\frac{%
\partial \mathbf{F}}{\partial R_{s}}\right) \rho +\frac{\partial
\mathbf{F}}{%
\partial \Omega }\omega \nonumber \\&& \hspace{6.5cm}=x\mathbf{y}^{\prime
\prime }+\mathbf{y}^{\prime }-%
\mathbf{F}.  \label{eq:10c} \\
&& \qquad \mathbf{z}^\prime(0)=-\mathbf{y}^\prime(0), \quad
\mathbf{z}(\infty)= -\mathbf{y}(\infty), \label{eq:10a}
\end{eqnarray}
\noindent where $\mathbf{E}$ is an identity $3\times 3$ matrix and
\begin{equation}
\>\dot{\mathbf{y}}=\mathbf{z}\>,\quad \dot{R}_{s}=\rho \>,\quad
\dot{\Omega}%
=\omega.  \label{eq:10b}
\end{equation}

The respective Frech\'{e}t derivatives at the point
$(\mathbf{y},R_{s},\Omega )$ are $\partial \mathbf{F}/ \partial
(.)$ and the dot in (\ref{eq:10b}) and below denotes the
differentiation with respect to ``time" $t$.

The solution $\mathbf{z}(x)$ of the above equation is sought as a
linear function towards the derivatives $\rho $ \ and $\omega $
\begin{equation}
\mathbf{z}=\mathbf{u}+\rho \mathbf{v}+\omega \mathbf{w},
\label{eq:11a}
\end{equation}
where $\mathbf{u}(x)\mathbf{,v}(x)$, and $\mathbf{w}(x)$ are
supposed to be new unknown 3D vector-functions of $x$.
Substituting for them in equation (\ref{eq:10c}), we obtain the
following three vector ODEs of second order with respect to these
quantities:
\begin{eqnarray}
-x\mathbf{u}^{\prime \prime }-\mathbf{u}^{\prime }+\frac{\partial
\mathbf{F}%
}{\partial \mathbf{y}}\mathbf{u}+\frac{\partial \mathbf{F}}{\partial
\mathbf{%
y}^{\prime }}\mathbf{u}^{\prime }\!\!\!&& =x\mathbf{y}^{\prime
\prime }+\mathbf{y}%
^{\prime }-\mathbf{F}        \label{eq:11} \\
-x\mathbf{v}^{\prime \prime }-\mathbf{v}^{\prime }+\frac{\partial
\mathbf{F}%
}{\partial \mathbf{y}}\mathbf{v}+\frac{\partial \mathbf{F}}{\partial
\mathbf{%
y}^{\prime }}\mathbf{v}^{\prime }\!\!\!&&
-\left( \frac{2}{R_{s}}\mathbf{F}+\frac{%
\partial \mathbf{F}}{\partial R_{s}}\right)  \label{eq:12} \\
-x\mathbf{w}^{\prime \prime }-\mathbf{w}^{\prime }+\frac{\partial
\mathbf{F}%
}{\partial \mathbf{y}}\mathbf{w}+\frac{\partial \mathbf{F}}{\partial
\mathbf{%
y}^{\prime }}\mathbf{w}^{\prime }\!\!\!&& =-\frac{\partial
\mathbf{F}}{\partial
\Omega }.        \label{eq:13}
\end{eqnarray}
The above three equations are coupled with the following six BCs:
\begin{eqnarray}
\mathbf{u}^{\prime }(0)=-\mathbf{y}^{\prime }(0)\!\!\!&& ,\qquad
\mathbf{u}(\infty
)=-\mathbf{y}(\infty );  \label{eq:14} \\
\mathbf{v}^{\prime }(0)=0\!\!\!&& ,\qquad \mathbf{v}(\infty )=0;
\label{eq:15} \\
\mathbf{w}^{\prime }(0)=0\!\!\!&& ,\qquad \mathbf{w}(\infty )=0,
\label{eq:16}
\end{eqnarray}
\noindent which are obtained from BCs (\ref{eq:10a}),
substituting for them with decomposition (\ref{eq:11a}) also. Let
us emphasize that the above equations (\ref{eq:11})-(\ref{eq:16})
have equivalent structures of the left-hand sides, which
essentially facilitates their numerical treatment.

In order to calculate the derivatives $\rho $ and $\omega $, we
apply CANM for the first additional BC (\ref{eq:4a}). This gives:
$$
{\dot{\sigma}}(0)=\sigma_c-\sigma (0).
$$
One more condition is required. Unfortunately, the second
additional condition (\ref{eq:4c}) is not convenient for this
purpose because knowledge about decomposition (\ref{eq:11a})
concerning the function $\mu (x)$ is not available. We avoid this
difficulty using the integral (\ref{eq:i1}) for $x_1 \equiv 0 $
and $x_2 \equiv 1$. Taking into account conditions (\ref{eq:8})
and (\ref{eq:4c}), we obtain an algebraic equation with respect
to the quantities $\nu(0)$, $\nu(1)$, $\varphi(0)$, $\varphi(1).$
After applying CANM to this equation, we get
\begin{eqnarray*}
&&{\dot{\nu}}(1)-{\dot{\nu}}(0)+2\,\frac{A^{\prime }\left[
\varphi (1)\right] }{A\left[ \varphi (1)\right]
}\,{\dot{\varphi}}(1)-2\frac{A^{\prime }\left[
\varphi (0)\right] }{A\left[ \varphi (0)\right] }\,{\dot{\varphi}}(0)= \\
&&\qquad\qquad\qquad \ln (1+\mu _c)-\left[ \nu (1)-\nu (0)\right]
-2\ln \frac{A\left[ \varphi (1)\right] }{A\left[ \varphi
(0)\right] }=0,
\end{eqnarray*}
\noindent where the abbreviation $A^\prime$ denotes the
derivative of the function $A$ with respect to the argument
$\varphi$.

Let us now eliminate all the derivatives in relation to ``time"
$t\,$\ by means of decomposition (\ref{eq:11a}). As a result, we
receive the following linear system of algebraic equations:
\begin{eqnarray}
&& a_{1}\rho +b_{1}\omega =c_{1} \nonumber\\ [-2mm]\label{eq:as}\\[-2mm]
&& a_{2}\rho +b_{2}\omega =c_{2} \nonumber
\end{eqnarray}
with respect to the unknown derivatives $\rho $\ and $\omega$. The
coefficients in formulae (\ref{eq:as}) are given by:
\begin{eqnarray*}
&& a_{1}=v_{1}(1)-v_{1}(0)+2\frac{A^{\prime }\left[ \varphi (1)\right] }{%
A\left[ \varphi (1)\right] } \, v_{2}(1)-2\frac{A^{\prime }\left[
\varphi
(0)\right] }{A\left[ \varphi (0)\right] } \, v_{2}(0) \\
&& b_{1}=w_{1}(1)-w_{1}(0)+2\frac{A^{\prime }\left[ \varphi (1)\right] }{%
A\left[ \varphi (1)\right] } \, w_{2}(1)-2\frac{A^{\prime }\left[
\varphi
(0)\right] }{A\left[ \varphi (0)\right] } \, w_{2}(0) \\
&& c_{1}=\ln \left[ 1+\mu_c \right] -\left[ \nu (1)-\nu (0)\right] -2\frac{%
A^{\prime }\left[ \varphi (1)\right] }{A\left[ \varphi (1)\right] }%
\,u_{2}(1)+2\frac{A^{\prime }\left[ \varphi (0)\right] }{A\left[
\varphi
(0)\right] }\,u_{2}(0) \\
&& \qquad \qquad -2\ln \frac{A\left[ \varphi (1)\right] }{A\left[ \varphi
(0)\right] }-u_{1}(1)+u_{1}(0) \\
&& a_{2}=v_{3}(0),\qquad b_{2}=w_{3}(0),\qquad
c_{2}=\sigma_c-\sigma (0)-u_{3}(0).
\end{eqnarray*}

Obviously, the explicit form of the coefficients in system
(\ref{eq:as}) depends on the concrete choice of functions
$f(\mu)$ and $g(\mu)$.

\section{General Sequence of the Algorithm}
We discretize the continuous ``time-like" parameter $t\in
[0,\infty )$ in the following way: $t_{k+1}=t_{k}+\tau _{k}$,
$t_{0}=0$, where $k=0,1,2,...$ denotes the number of iterations,
and the ``time" step $\tau _{k}$ is generally assumed as a
variable quantity. Next, we use the Euler difference scheme
\cite{jmp} to approximate the ``time" derivatives in equations
(\ref{eq:10b}). Then we can write:
\begin{eqnarray}
\mathbf{y}_{k+1}(x)\!\!\!&&=\mathbf{y}_{k}(x)+\tau _{k}\left[ \mathbf{u}%
_{k}(x)+\rho _{k}\mathbf{v}_{k}(x)+\omega _{k}\mathbf{w}_{k}(x)\right] ,
\nonumber \\
R_{s,k+1}\!\!\!&&=R_{s,k}+\tau _{k}\rho _{k},        \label{eq:22} \\
\Omega _{k+1}\!\!\!&&=\Omega _{k}+\tau _{k}\omega _{k}.  \nonumber
\end{eqnarray}
Let us suppose that the functions $\nu _k (x)$, $\varphi _k(x)$,
$\sigma_k(x) $, $\mu _k(x)$ and the parameters $R_{s,k}$, $\Omega
_{k}$ are given. We solve the linear BVP
(\ref{eq:11})-(\ref{eq:13}) and, thus, we compute the functions
$\mathbf{u}_{k}(x)$, $\mathbf{v}_{k}(x)$, $\mathbf{w}_{k}(x)$.
Next, to obtain the derivatives $\rho _{k}$\ and $\omega _{k}$\
we solve system (\ref{eq:as}). After that, using decomposition
(\ref{eq:22}) for a selected $\tau_k$, we calculate the functions
$\nu _{k+1}(x)$, $\varphi _{k+1}(x)$, $\sigma_{k+1}(x)$, the
radius of the star $R_{s,k+1}$, and the quantity $\Omega _{k+1}$
as well at the new stage $k+1$. In the end, we calculate the
function $\mu_{k+1}(x)$ at the new stage, according to the
recurrent formula, which can be obtained immediately from the
first integral (\ref{eq:i1}).

For every iteration $k$\, an optimal time step $\tau _{opt}$ is
determined in accordance to the Kalitkin\& Ermakov formula
\cite{ermakov}, \cite{zhanlav}:
\begin{equation}\label{ek}
  \tau _{opt}=\frac{\delta (0)}{\delta (0)+\delta (1)},
\end{equation}
where the residual $\delta (\tau )$ is represented as follows:
\[
\delta (\tau _{k})=\max \left[ \delta _{f},(R_{s,k}+\tau _{k}\rho
_{k})^{2},(\Omega _{k}+\tau _{k}\omega _{k})^{2}\right]
\]
and $\delta _{f}$ is the Euclidean residual of RHS of the
equation (\ref{eq:11}). Formula (\ref{ek}) provides approximately
the minimal value of the residual for the current solution, given
by (\ref{eq:22}).

The criterion for termination of the iterations is $\delta (\tau
_{opt})<\varepsilon $, where $\varepsilon \sim 10^{-8}\div
10^{-12}$. Then, for the sought solutions we set $\nu(x) \equiv
\nu_{k+1}(x)$, $\varphi(x) \equiv \varphi_{k+1}(x)$, $\sigma(x)
\equiv \sigma_{k+1}(x)$, $R_s \equiv R_{s,k+1}$, $\Omega \equiv
\Omega_{k+1}$.

The use of the standard programs available, for example, via the
Internet \cite{nanet}, to solve numerically the linear BVPs
(\ref{eq:11})-(\ref{eq:16}) is unhandy for many reasons. Because
of that, the spline-collocation scheme is employed in our case.

We introduce a nonuniform grid
$$\Delta: \> x_{i+1}=x_{i}+h_{i},\quad i=0,1,\ldots ,N_{s},N_{s+1},\ldots,
N-1,
\quad x_0=0, \quad x_N=X_{\infty },$$ \noindent on the interval
$x\in [0,X_{\infty }]$, condensing to the points $x=0$ and $x =
1$. Here, $X_{\infty}$ is the ``actual infinity", $N_{s}$ is the
number of the node $x=1$, $N$ is the full number of the
subintervals, and $h_{i}$ is the grid step. We will seek
approximate solutions of the above linear BVPs as a cubic spline
on the grid $\Delta$. Namely, for $x\in [x_{i},x_{i+1}]$, $\>
i=0,...,N-1$ we set
\begin{equation}\label{spl}
\mathbf{U}(x)=\psi_{1}(\theta )\,\mathbf{U}_{i}+\psi_{2}(\theta )\,\mathbf{M
}
_{i}+\psi_{3}(\theta )\,\mathbf{U}_{i+1}+\psi_{4}(\theta
)\,\mathbf{M}_{i+1}.
\end{equation}
\noindent In the above formula the relative coordinate $\theta=(x
- x_i)/h_i$ and the known functions $\psi_l\>(\theta ),\> l =
1,\dots, 4$, are the coefficients of the spline. For simplicity
in the last formula, we introduced the $3\times 3$ matrices $
\mathbf{U}$ and $ \mathbf{M}$, consisting of the coordinates of
the vectors $\mathbf{u},\mathbf{v},\mathbf{w}$ from
(\ref{eq:11a}) and their first moments at the spline nodes $x_i
\>$, $i=0,\ldots,N$. According to the collocation method
\cite{DeBoor}, in every subinterval $[x_{i},x_{i+1}]\>,$
$i=0,...,N-1,$ the system (\ref{eq:11}) - (\ref{eq:13}) is
satisfied at the corresponding Gaussian points
$\theta_1=1/2-\sqrt{3}/6$ and $\theta_2=1/2+\sqrt{3}/6$. This
kind of discretization yields an algebraic system with respect to
the functions and their moments at the spline nodes. The
corresponding matrix has an almost block-diagonal structure (see
\cite{DeBoor}). Therefore, at the \textit{i}$-$th block
($i=1,\ldots ,N-1$) the collocation equations have the form:
\[
\left(\begin{array}{cccc} \Vert a_{kn}^1 \Vert & \Vert b_{kn}^1
\Vert & \Vert c_{kn}^1 \Vert & \Vert
d_{kn}^1 \Vert \\
&&& \\
\Vert a_{kn}^2 \Vert & \Vert b_{kn}^2 \Vert & \Vert c_{kn}^2 \Vert
& \Vert d_{kn}^1 \Vert
\end{array} \right)
\left( \begin{array}c
{\bf U}_i \\
{\bf M}_i \\
{\bf U}_{i+1} \\
{\bf M}_{i+1}
\end{array}\right) =%
\left(\begin{array}c
{\bf e}_i^1 \\
\\
{\bf e}_i^2
\end{array} \right),
\]
\noindent where $\mathbf{e}_{i}$ is the vector of RHSs of the
equations (\ref{eq:11})-(\ref{eq:13}) at the collocation nodes,
while the superscript corresponds to the number of these nodes.
Here:
\begin{eqnarray*}
&&
_{kn}^j=-\left( \frac{\xi_{ij}}{h_{i}^{2}}{\ddot{\psi}}_{1j}+\frac{1}{h_{i}}
{%
\dot{\psi}}_{1j}\right) \delta _{kn}+\left( \frac{\partial
F_{k}}{\partial
y_{n}^{\prime }}\right) _{ij}\frac{1}{h_{i}}{\dot{\psi}}_{1j}+\left( \frac{%
\partial F_{k}}{\partial y_{n}}\right) _{ij}\psi_{1j}, \\
&& b_{kn}^j=-\left(
\frac{\xi_{ij}}{h_{i}}{\ddot{\psi}}_{2j}+{\dot{\psi}}_{2j}\right)
\delta _{kn}+\left( \frac{\partial F_{k}}{\partial y_{n}^{\prime
}}\right) _{ij}{\dot{\psi}}_{2j}+\left( \frac{\partial
F_{k}}{\partial y_{n}}\right)
_{ij}h_{i}\psi_{2j}, \\
&&
_{kn}^j=-\left( \frac{\xi_{ij}}{h_{i}^{2}}{\ddot{\psi}}_{3j}+\frac{1}{h_{i}}
{%
\dot{\psi}}_{3j}\right) \delta _{kn}+\left( \frac{\partial
F_{k}}{\partial
y_{n}^{\prime }}\right) _{ij}\frac{1}{h_{i}}{\dot{\psi}}_{3j}+\left( \frac{%
\partial F_{k}}{\partial y_{n}}\right) _{ij}\psi_{3j}, \\
&& d_{kn}^j=-\left(
\frac{\xi_{ij}}{h_{i}}{\ddot{\psi}}_{4j}+{\dot{\psi}}_{4j}\right)
\delta _{kn}+\left( \frac{\partial F_{k}}{\partial y_{n}^{\prime
}}\right) _{ij}{\dot{\psi}}_{4j}+\left( \frac{\partial
F_{k}}{\partial y_{n}}\right)
_{ij}h_{i}\psi_{4j}, \\
&&\qquad \hbox {for}\qquad  k=1,2,3,\quad n=1,2,3,\quad j=1,2,
\end{eqnarray*}
\noindent and the quantities $\xi_{ij}=x_i+\theta_j h_i$ are the
absolute coordinates of the collocation points. The derivatives of
the spline coefficients $\psi_l$ with respect to the relative
coordinate $\theta$ are dotted.

The dimensions of the first and the last blocks in the global
matrix are greater, since we add two matrix rows corresponding,
respectively, to the left and right BCs.

Formula (\ref{spl}) is also used for the approximation of the
RHSs of system (\ref{eq:11})~-~(\ref{eq:13}) in the collocation
points.

The spline-difference schemes of this kind have a high order of
approximation $\mathcal{O}({\bar{h}^{4}})\ $, where
${\bar{h}}=\max \{h_{i}\},$ $i=0,...,N$.

It is clear that for solving all the three algebraic systems,
corresponding to the linear BVPs (\ref{eq:11}) - (\ref{eq:16}) at
every iteration, only one $LU$-decomposition is necessary.

Depending on the initial values of the government physical
parameters, the number of iterations varies approximately in the
range 4 $\div $ 16. If we vary some solution as a function of one
of the parameters $\mu_c$, $\sigma_c$, $\gamma $, $\Lambda$, or
$b$, then we use the previous solution as an initial
approximation for computing the next one.

\section{Results and Discussion}
In order to be specific in the present article, we focus our
attention on a concrete scalar-tensor gravity model, characterized
by the functions
\[
A(\varphi )=\exp (\frac{\varphi }{\sqrt{3}})\quad \hbox{and}\quad
V(\varphi )=(1-[A(\varphi )]^{-1})^2.
\]
For more details concerning this gravitational model, we refer the
reader to the recent paper \cite{F} and the references therein.

The order of approximation of the used spline-difference scheme is
verified by the Runge rule.

The Runge rule is presented by the formula:
\[
\frac{y_{h}-y_{\frac{h}{2}}}{y_{\frac{h}{2}}-%
y_{\frac{h}{4}}} = 2^{p},
\]
where $p$ is Runge's number and $y_h$, $y_{\frac{h}{2}}$, $
y_{\frac{h}{4}}$ are the values of the grid function $y$ at the
given node, computed on meshes with steps $h, h/2$, and $h/4$. In
our case $p$ must be approximately equal to 4.

In Table 1 the values of the sought grid functions at the point
$x=1$, the corresponding radius of star $R_{s}$, and the quantity
$\Omega $ for $\sigma_c=0.8$, $\mu_c=1$, $\Lambda =0.01$, $\gamma
=1$, $b=1$, and $X_{\infty} = 128$, are shown.

\begin{table}[ht]
\caption{ Data for checking the Runge rule.}
%% way to get table to spread out to width of page:
\begin{tabular*}{\textwidth}{@{\extracolsep{\fill}}cccccc}
\hline $h$ & $\nu (1)$ & $\varphi (1)$ & $\sigma (1)$ & $R_{s}$ &
$\Omega $\cr \hline
$\frac{1}{16}$ & $-1.0059230404$ & $-0.0471137759$ & $0.4777335163$ & $%
1.1609111685$ & $0.8006662485$ \cr
$\frac{1}{32}$ & $-1.0059334054$ & $-0.0471120738$ & $0.4777483180$ & $%
1.1608888836$ & $0.8006671950$ \cr
$\frac{1}{64}$ & $-1.0059342032$ & $-0.0471119781$ & $0.4777490917$ & $%
1.1608875328$ & $0.8006672467$ \cr \hline $p$ & $3.61$ & $4.22$ &
$4.37$ & $4.06$ & $4.28$ \cr \hline
\end{tabular*}
\end{table}

Therefore, it is obvious that the Runge relationship is satisfied
both for the functions and the eigenvalues $R_{s}$\ and $\Omega $.

The correctness of the spline-difference scheme is verified
through appropriate numerical experiments consisting of both grid
doubling and doubling of the ``actual infinity". For this
purpose, uniform meshes are used with numbers of the spline nodes
$N= 256$, $512$, $1024$, $2048$, respectively. It turns out that
the relative error between the values of the functions $\nu (x)$,
$\varphi (x)$, and $\sigma (x)$, varies in the range $0.1\%-1\%$
when the mesh is ``coarse" ($N=256, 512$), and in the range
$0.003\%-0.02\%$ when the mesh is ``fine" ($N=1024,$ $2048$).
Similar experiments are carried out with the ``actual infinity"
$X_{\infty }=64,128,256$. It is interesting to note that the
relative error between the set functions $\varphi (x)$ and $\sigma
(x)$ is very small (less than $10^{-4}\,\%$), while the function
$\nu (x)$ is more sensitive with respect to the choice of the
quantity $X_{\infty }$. This fact is fully explainable if we take
into account that the function $\nu (x)$ decreases slowly at the
infinity compared to the other functions. (Theoretically $\nu
(x)\sim -\frac{M}{R_s x}$ when $x\to \infty $. Here, the quantity
$M$ is the total star mass.) The computed values of the
derivative $\nu ^{\prime }(X_{\infty })$ as a function of the
``actual infinity" $X_{\infty }$ are presented in Table 2. It is
easy to see the relationship $\nu ^{\prime }(X_{\infty
})=\frac{C}{X_{\infty }^{2}}$, where the constant $C > 0$ depends
on the concrete solution (for the above solution $%
C\approx 1.133$).

\begin{table}[ht]
\caption{ Asymptotic behaviour of the derivative $\nu^\prime$ at
the ``actual infinity" $X_{\infty}$.}
\begin{tabular*}{\textwidth}{@{\extracolsep{\fill}}cccccc}
\hline $X_{\infty }$ & $32$ & $64$ & $128$ & $256$ & $512$\cr
\hline $\nu ^{\prime }(X_{\infty })$ & $1.07246 \!\!\times\!\!
10^{-3}$ & $2.63721 \!\!\times\!\! 10^{-4}$ & $6.53945
\!\!\times\!\! 10^{-5}$ & $1.62825 \!\!\times\!\! 10^{-5}$ &
$4.06241 \!\!\times\!\! 10^{-6}$ \cr \hline
\end{tabular*}
\end{table}

All government parameters are varied in wide
physically-admissible ranges. As initial distributions of the
functions $\nu (x)$, $\varphi (x)$, $ \sigma (x)$ and $\mu (x)$
both analytic and numerical approximations are used.

Results concerning a family of solutions will be considered below.
They are obtained for the following fixed values of the
parameters: $\mu _c=0.5$, $ \Lambda =10$, $\gamma =10$, $b=1$,
and the ``actual infinity" $X_{\infty }=128$, when the parameter
$\sigma_c$ runs the interval $[0.1,0.9]$.

\begin{figure}[ht]
\centerline{\psfig{figure=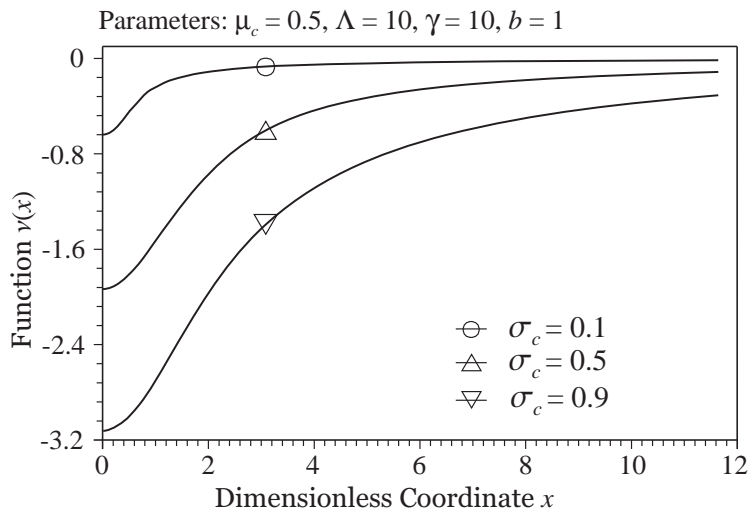,width=4.in}} \vskip12pt
\caption
{The function $\nu (x)$ in dependence on the parameter $%
\sigma_c$: ``$\circ $" - $\sigma_c=0.1;$ ``$\triangle $" -
$\sigma_c=0.5;$ ``$\nabla $" - $\sigma_c=0.9.$}
\end{figure}

Figure 1 presents the dependence of the function $\nu (x)$ on the
dimensionless coordinate $x$ for three different values of the
central bosonic density $\sigma _c$. It is seen that when
$\sigma_c$ increases, the absolute value of $\nu (x)$ as a whole
decreases and at great distances (from 3 star radii when
$\sigma_c=0.1$ until 45 star radii in the case $\sigma_c=0.9$)
from the star's center approaches asymptotically zero. The
qualitative behaviour of the three curves, however, remains the
same. Such a behaviour is natural and should be expected if the
differential equation (\ref {eq:1}) for $\nu (r)$ is taken into
account. From physical point of view, this behaviour is natural
also because the function $\exp(\frac{\nu (x)}{2})$ is related to
the gravitational potential.

\begin{figure}[ht]
\centerline{\psfig{figure=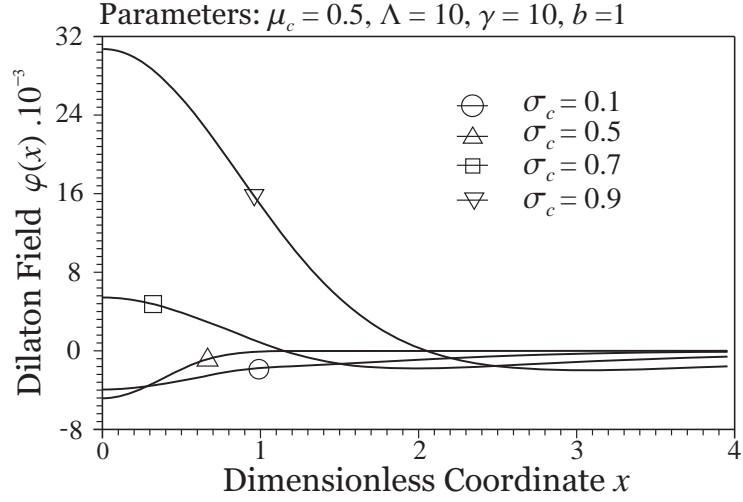,width=4.in}} \vskip12pt
\caption {The potential of dilaton $\varphi (x)$ as function of
the parameter $\sigma_c$: ``$\circ $" - $\sigma_c=0.1;$
``$\triangle $" - $\sigma_c=0.5;$ ``$\Box $" - $\sigma_c=0.7;$
``$\nabla $" - $\sigma_c=0.9$.}
\end{figure}

Figure 2 presents the dependence of the dilaton field $\varphi
(x)$ on the dimensionless coordinate $x$ for four different values
of $\sigma_c$. The qualitative behaviour of the field $\varphi
(x)$ as a function of $\sigma_c$ is the following. For small
values when $\sigma_c$ increases, the dilaton field around the
center of the star decreases. Then, after some critical value
$\sigma_c^{*}$ the behaviour of $\varphi (x)$ is changed and
$\varphi (x)$ around the center of the star begins to increase
with the increase of $\sigma_c$. The cause of the described
behaviour is the presence of the term ${\stackrel{B}{T}}$ on the
RHS of equation (\ref{eq:2}). For sufficiently small values of
the density $\sigma_c$ the term ${\stackrel{B}{T}}$ is negative
and has a dominant contribution with respect to the term
${\stackrel{F}{T}}$. For the sufficiently large central value
$\sigma _c$ ($\sigma_c\ge \sigma_c^{*}$), the term ${\
 \stackrel{B}{T}}$ changes its sign and amplifies the contribution of ${
\stackrel{F}{T}}$, leading to the increase of the function
$\varphi (x)$.

\begin{figure}[ht]
\centerline{\psfig{figure=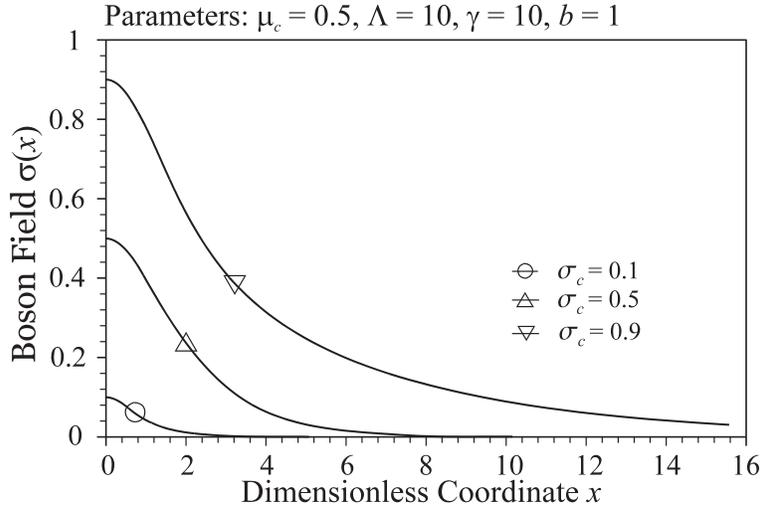,width=4.in}} \vskip12pt
\caption
{The bosonic density $\sigma (x)$ in dependence on the parameter $%
\sigma_c$: ``$\circ $" - $\sigma_c=0.1;$ ``$\triangle $" -
$\sigma_c=0.5;$ ``$\nabla $" - $\sigma_c=0.9.$}
\end{figure}

From a physical point of view, the described behaviour of the
dilaton field (and consequently  the behaviour of the physical
gravitational ``constant" $G_{*}A^{2}(\varphi)\>$) for the
central values $\sigma_c> \sigma_c^{*}$ seems to be strange. In
order to clarify this situation, we have to take into account that
in the range $\sigma_c> \sigma_c^{*}$ (for the fixed value of the
central fermionic density $\mu_c$) the star is unstable and,
therefore, the mentioned range is not physically relevant. Such a
behaviour has to be considered only as an iteresting mathematical
fact. In the domain of stability $0 < \sigma_c< \sigma_c^{*}$, as
we have already seen, the dilaton field $\varphi (x)$ has a normal
physical behaviour - it decreases when the parameter $\sigma_c$
increases.

\begin{figure}[ht]
\centerline{\psfig{figure=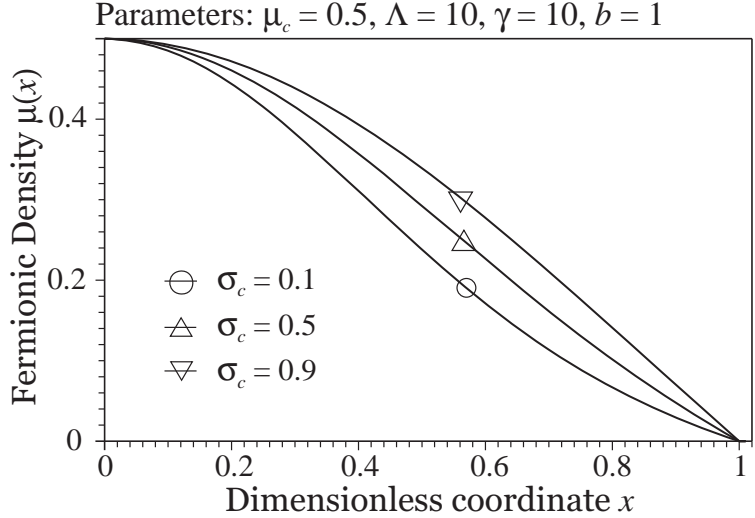,width=4.in}} \vskip12pt
\caption
{The fermionic density $\mu (x)$ in dependence on the parameter $%
\sigma_c$: ``$\circ $" - $\sigma_c=0.1;$ ``$\triangle $" -
$\sigma_c=0.5;$ ``$\nabla $" - $\sigma_c=0.9.$}
\end{figure}

The dependence of the bosonic density $\sigma (x)$ on the
dimensionless coordinate $x$ for three different values of
$\sigma_c$ is presented on Figure 3. The qualitative behaviour is
the same for all three different values of $\sigma_c$. It
approaches zero at infinity (rapidly when $\sigma_c=0.1$ and more
slowly when $\sigma_c$ increases).

\begin{figure}[ht]
\centerline{\psfig{figure=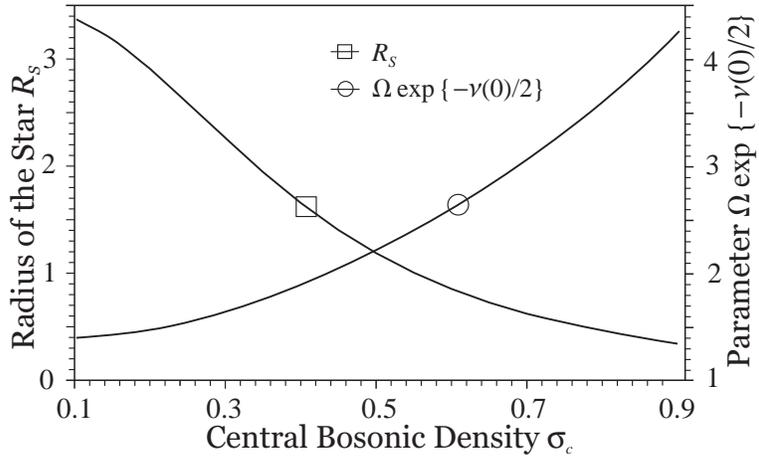,width=4.in}} \vskip12pt
\caption {The radius of the star $R_{s}$ and the quantity $\Omega
\exp (-\dfrac{\nu (0)}{2})$ as functions of the parameter $\sigma_c$: $%
R_{s}$ - ``$\Box $"; $\Omega\exp (-\dfrac{\nu (0)}{2})$ -
``$\circ$".}
\end{figure}

In the next Figure 4 the dependence of the fermionic density
$\mu(x) $ on the dimensionless coordinate $x$ is presented for
three different values of $\sigma_c$. The qualitative behaviour
of the three curves is similar. In agreement with the initial
assumption, it is nontrivial only within the star. It is seen
that when the value of $\sigma_c$ increases, the density $\mu (x)$
increases as a whole, too. This fact is related to the effect of
an increase of the gravitational field with the increase of
$\sigma_c$ - the star becomes more compact, which leads to the
greater density of matter, respectively to the function $\mu(x) $.
The same may be seen in Figure 5 - when the central value
$\sigma_c$ increases, the radius of the star $R_{s}$ decreases
about 10 times.

From physical point of view it is important to get knowledge about
the behaviour of the quantity $\Omega \exp (-\frac{\nu (0)}{2})$
as a function of the central value $\sigma_c$. That quantity may
be considered as the energy of one boson particle in the
gravitational field yielded by the rest matter (in the Einstein
frame). Figure 5 clearly shows that the quantity $\Omega \exp
(-\frac{\nu (0)}{2})$ increases along with $\sigma_c$. Such
behaviour should be expected, because the energy of the system has
to increase with the increase of the central density $\sigma_c$
of the star.

\section*{Concluding Remarks}

Based on CANM an iteration method for solving the nonlinear BVP,
describing a static spherically-symmetric boson-fermion star, is
developed.

A linearization of the main equations of the star renders the
original two-parametric nonlinear spectral problem to three
two-point linear vector BVPs and a linear system of algebraic
equations for the spectral parameters (the radius of the star
$R_s$ and the frequency $\Omega$ of the bosonic field). A
spline-collocation scheme of fourth order of approximation for
solving numerically these BVPs is used.

Our basic physical result is that the structure and the
properties of the star in the presence of a massive dilaton field
depend essentially on both its fermionic and bosonic components.
This shows that a careful investigation of these properties may
give new real ways to discover physical effects of the
hypothetical boson fields and dilaton field in stars.

\section*{Appendix}

For reader's convenience, we briefly explain the main ideas of
CANM.

CANM can be treated as a particular case of the continuous
analogues of iteration methods, strictly formulated and studied
by M.K. Gavurin in 1958 (see the review in \cite{gavurin}). Among
the number of papers devoted to the theoretical development and
applications of CANM for solving wide classes of nonlinear
equations, we will indicate the basic papers \cite{pon72} -
\cite{zhanlav} as well as the reviews \cite{jmp}, \cite{puz99}.

Let us consider the nonlinear equation:
\begin{equation}
\label{a1} \chi (y) = 0,
\end{equation}
\noindent where $\chi(y)$ is an operator defined in a Banach space
$\mathbf{Y}$. We suppose that the equation (\ref{a1}) has an
isolated exact solution $y^{*} \in \mathbf{Y}$. Let the element
$y_0\in \textbf{Y}$ (an initial approximation to $y^{*}$) be
given. To solve equation (\ref{a1}), we can use an iteration
process, usually taking it in the form:
$$  y_{n+1} = y_n + \psi(y_n), \quad n = 0,1, 2, ....$$
\noindent Here, $n$ indicates the number of iterations and $\psi$
is an appropriate function, which carries $\mathbf{Y}$ into itself
and has the same zeroes as $\chi$.

The choice of the function $\psi(y)$ depends on the kind of
concrete iteration method used.

According to Gavurin's idea, for each iteration process of such
kind it is possible to formulate the corresponding continuous
analogue in the following way. Let us consider an abstract
function $y(t)$ of the independent continuous variable
$t\in[0,\infty)$ instead of the sequence $y_0, y_1, ..., y_n,
...$, and suppose that $y(t_n)=y_n$ for each $n$. Then, we can
introduce the derivative $\dot{y}(t)$ instead of the increment
$y_{n+1}-y_n$ and replace (\ref{a1}) with the abstract initial
value problem on the interval $t \in [0,\infty)$
\begin{equation} \label{a3}
\dot{y}(t) = \psi(y), \quad y(0) = y_0.
\end{equation}

Such a transition from a difference equation to a differential one
has many advantages, both in pure theoretical and applied aspects.

In the case of Newton's method, we set $\psi(y) = -
\chi^\prime(y)^{-1} \chi(y)$, where $\chi^\prime(y)$ is the
corresponding Frech\'{e}t derivative of $\chi(y)$. Then, the main
equation of CANM, arising from (\ref{a3}), can be rewritten in
the form:
\begin{equation}\label{a4}
  \chi^\prime (y) \dot{y} = -\chi(y).
\end{equation}

Obviously, the above ODE has a significant first integral of the
kind:
\begin{equation}
\label{a5} \chi \left( {y\left( {t} \right)} \right) = \chi \left(
{y_{0}}  \right)e^{ - t},
\end{equation}
\noindent which means that $\chi(y(t)) \to 0$ when $ t \to \infty
$.

Various theorems, based on (\ref{a5}), concerning the convergence
of a path $y(t)$ to the exact solution $y^*$ have been proved. For
example, a theorem \cite{jmp}, which guarantees the convergence of
CANM for a simple BVP, is cited below.

The following BVP is considered:
\begin{eqnarray}
&&  -y^{\prime \prime}+f(x,y)=0, \quad x\in(0,1), \label{a8}\\
&& \qquad y(0)= 0, \quad y(1)=0. \label{a9}
\end{eqnarray}
\begin{theorem}
Let the BVP (\ref{a8}), (\ref{a9}) have an isolated solution
$y^*(x)$ and:

\begin{description}
\item[i) ] the function $f(x,y)$ have continuous partial derivatives up to
the second order in some domain $D$;

\item[ii) ] the linear BVP
\begin{eqnarray*}
&& -w^{\prime \prime}+f_y^\prime(x,y) w = 0, \quad x\in (0,1),\\
&& \qquad w(0) = 0, \quad w(1) =0,
\end{eqnarray*}
have only a trivial solution $w(x)\equiv 0$ for every smooth
function $y(x)\in D$;

\item[iii) ] the initial approximation $y_0(x)\in D$ be a smooth enough
function satisfying:
$$\|-y_0^{\prime\prime}+f(x,y_0\|\leq\varepsilon
\quad \hbox{for} \quad \varepsilon > 0.$$
\end{description}

Then the system
$$ -w^{\prime \prime}+f_y^\prime(x,y) w = y^{\prime
\prime}-f(x,y), \quad \dot{y}=w,$$ with BCs $\> w(0,t)=0,\,
w(1,t)=0,$ and an initial condition $y(x,0)=y_0(x),$ has in
$[0,1]\cup [0,\infty)$ an unique solution, satisfying the
relation:
$$\lim_{t \to\infty} \|y(x,t)-y^{*}(x)\|_{C^2[0,1]}=0.$$
\end{theorem}
The numerical solution of CANM equation (\ref{a4}) is based on an
appropriate scheme for discretization, which has to be stable for
the asymptotic stability of the path $y(t)$. The most frequently
used one is Euler's scheme (see the details in the above cited
papers). At first, the linearized equation:
\begin{equation}\label{a6}
\chi^\prime (y_n) w_n = -\chi(y_n),
\end{equation}
\noindent is solved with respect to the increment $w_n$, and then
the next approximation is obtained via the formula:
\begin{equation}\label{a7}
y_{n+1} = y_n + \tau_n w_n.
\end{equation}
Here, $0 < \tau_n \leq 1$ is an iteration parameter. When $\tau_n
\equiv 1$, the classical Newton method is obtained. We note that
the choice of $\tau_n$ is important for the rapid convergence of
the process. It is possible to choose this parameter so that the
range of convergence is wider in comparison to the classical
Newton's method \cite{ermakov}, \cite{zhanlav}.

Theorems regarding the convergence of iterations (\ref{a6}),
(\ref{a7}) for wide enough hypotheses as well as essential
generalizations of CANM, are discussed in the above cited papers.

\bigskip
{\bf {Acknowledgment.}}
We thank Prof. Igor V. Puzynin (JINR, Dubna, Russia) for useful
remarks.

\end{document}